\documentclass[reqno,11pt]{article}
\usepackage{amsmath,amssymb}
\usepackage{amsthm,amscd,amsfonts,}
\usepackage{amsmath}
\usepackage{amsthm}
\usepackage{graphicx}
\usepackage{amsmath,amsthm}
\usepackage{amsmath,amstext,amsfonts}
\usepackage{subfigure}
\usepackage[T1]{fontenc}
\usepackage[utf8]{inputenc}
\usepackage{authblk}
\newtheorem{theorem}{Theorem}[section]

\newtheorem{corollary}{Corollary}[section]
\newtheorem{example}{Example}[section]

\newcommand{\R}{\mathbb{R}}
\numberwithin{equation}{section}

\begin{document}

\title{\textbf{Integrating Factors and First Integrals for a Class of Higher Order Differential Equations}}
\author{M. Al-Jararha\thanks{mohammad.ja@yu.edu.jo}}
\affil{Department of Mathematics, Yarmouk University, Irbid, Jordan, 21163.}

\date{}
\maketitle

\vspace{0.2in}

\noindent \textbf{AMS (2000) Subject Classification}: 34A25, 34A30.

\begin{abstract}
If the $n-th$ order differential equation is not exact, under certain conditions, an integrating factor exists which  transforms the differential equation into an exact one. Hence, its order can be reduced to the lower order. In this paper, the principle of finding an integrating factor of a none exact differential equations is extended to the class of $n$-th order differential equations
\begin{align}
F_n\left(t,y,y^\prime,y^{\prime\prime},\ldots,y^{(n-1)}\right)y^{(n)}&+F_{n-1}\left(t,y,y^\prime,y^{\prime\prime},\ldots,y^{(n-1)}\right)y^{(n-1)}+\cdots +\nonumber\\
&+F_{1}\left(t,y,y^\prime,y^{\prime\prime},\ldots,y^{(n-1)}\right)y^{\prime}+F_{0}\left(t,y,y^\prime,y^{\prime\prime}\ldots,y^{(n-1)}\right)\nonumber\\
&=0,\nonumber
\end{align}
where $F_0,F_1,F_2, \cdots,F_n$ are continuous functions with their first partial derivatives
on some simply connected domain $\Omega \subset\R^{n+1}$.  In particular, we prove some explicit forms of integrating factors for this class of differential equations. Moreover, as a special case of this class, we consider the class of third order differential equations in more details. We also present some illustrative examples.
\end{abstract}
\maketitle
%---------------------------------------------------------------------
\section{Introduction}

Differential equations play a major role in Applied Mathematics, Physics, and Engineering \cite{Ames, Jordan, Lefschetz, Struble,von}. To find the general solution of a differential equation is not an easy problem in the general case. In fact, a very specific class of differential equations can be solved by using special techniques and transformations. One of these techniques is to reduce the order of the differential equation by finding a proper integrating factor. Recently, many studies appear to deal with the problems of the existence and finding integrating factors of certain classes of differential equations. For example, in \cite{AlAhmad, existencepp, Cheb, Bouquet}, the authors investigated the existence of integrating factors for some classes of second order differential equations. In \cite{existencepp}, the authors investigated the existence of  integrating factors of $n$-th order system differential equations which has known symmetries of certain type. 

In \cite{kai}, the authors improve some symbolic algorithms to compute the integrating factor for certain class of third order  differential equations. In this paper, we investigate the existence of integrating factors of the following class of $n-th$ order  differential equations: 
\begin{align}\label{HODE}
F_n\left(t,y,y^\prime,y^{\prime\prime},\ldots,y^{(n-1)}\right)y^{(n)}&+F_{n-1}\left(t,y,y^\prime,y^{\prime\prime},\ldots,y^{(n-1)}\right)y^{(n-1)}+\cdots+ \nonumber\\
&+F_{1}\left(t,y,y^\prime,y^{\prime\prime},\ldots,y^{(n-1)}\right)y^{\prime}+F_{0}\left(t,y,y^\prime,y^{\prime\prime}\ldots,y^{(n-1)}\right)\nonumber\\
&=0.
\end{align}
where $F_0,F_1,F_2, \cdots,F_n$ are continuous functions with their first partial derivatives 
on some simply connected domain $\Omega \subset\R^{n+1}$. In fact, we prove some theoretical results related to the existence of certain forms of integrating factor for the class of differential equations  \eqref{HODE} as well as we give explicit forms for such integrating factor. We also present some illustrative examples.  

The paper layout: In section \textbf 2, we prove the main result. In section \textbf 3 is devoted for concluding remarks. 
%---------------------------------------------------------------------
\section{Integrating Factor and First Integrals for a Class of $n$-th Order Differential Equations}
In this section, we investigate the existence of integrating factors of
 \eqref{HODE} when it is not an exact equation. Generally, the $n$-th order differential equation $$f(t,y,y^\prime,\cdots,y^{(n-1)},y^{(n)})=0$$ is called  exact if there exists a differentiable function $\Psi(t,y,y^\prime,\cdots,y^{(n-1)})=c$, such that $\frac{d}{dt}\Psi(t,y,y^\prime,\cdots,y^{(n-1)})=f(t,y,y^\prime,\cdots,y^{(n-1)},y^{(n)})=0$. In this case, $\Psi(t,y,y^\prime,\cdots, y^{(n-1)})=c$ is called the first integral of  $f(t,y,y^\prime,\cdots,y^{(n-1)},y^{(n)})=0$, e.g., see,   \cite{cov,Peter}. In \cite{Aljararha}, the author gave the explicit conditions for 
\eqref{HODE}
to be exact. He also gave an explicit formula for the first integral $\Psi\left(t,y,y^\prime,\cdots,y^{(n)}\right)=c$. Particularly, we have the following theorem:
%%%%%%%%%%%%%%%%%%%%%%%%%%%THM%%%%%%%%%%%%%%%%%%%%%%%%%%%%%%%%%%%%%%%
\begin{theorem}\emph{\cite{Aljararha}}\textbf .
Assume that $F_0,F_1,F_2,\ldots,F_n$ are continuous with their first partial derivatives on a simply connected domain $\Omega$ in $\R^{n+1}$. Then the differential equation \eqref{HODE} is exact if 
 $\displaystyle \frac{\partial F_i}{\partial y^{(j-1)}}=\frac{\partial F_j}{\partial y^{(i-1)}}$, for all $i=2,3,\ldots,n$ and $j=1,2,\ldots,i-1$
 and $\displaystyle \frac{\partial F_i}{\partial t}=\frac{\partial F_0}{\partial y^{(i-1)}}$, for all $i=1,2,3,\ldots,n$. Moreover, its first integral is explicitly given by 
 \begin{eqnarray*}%\label{psinthdorder}
 \Psi\left(t,y,y^\prime,\cdots,y^{(n-1)}\right)&=& \int_{t_0}^{t}F_0\left(\eta,y,y^\prime,\cdots,y^{(n-1)}\right)d\eta+  \int_{y_0}^{y}F_1\left(t_0,\eta,y^\prime,\cdots,y^{(n-1)}\right)d\eta\nonumber\\
&+&\cdots+\int_{y_0^{(n-1)}}^{y^{(n-1)}}F_n\left(t_0,y_0,y^\prime_0,\cdots,\eta\right)d\eta=c,
\end{eqnarray*}
where  $c$ is an integrating constant. $\square$
\end{theorem}
%%%%%%%%%%%%%%%%%%%%%%%%END THM%%%%%%%%%%%%%%%%%%%%%%%%%%%%%%%%%%%
 Assume that \eqref{HODE} is a none exact differential equation. Then according to the above theorem, an integrating factor $\mu(t,y,y^\prime,\cdots,y^{(n-1)})$ of \eqref{HODE}  exists if it solves the following system of $n!$ first order partial differential equations:
\begin{equation}\label{sys1}
\left\{
\begin{array}{ll}
\displaystyle \mu(\mathbf y)\frac{\partial F_{i}(\mathbf y)}{\partial t}+\frac{\partial \mu(\mathbf y)}{\partial t}F_{i}(\mathbf y )=\mu(\mathbf y )\frac{\partial F_0(\mathbf y)}{\partial y^{(i-1)}}+\frac{\partial \mu(\mathbf y )}{\partial y^{(i-1)}}F_{0}(\mathbf y),\; i=1,2,\cdots,n,\\\\
\displaystyle \mu(\mathbf y)\frac{\partial F_{i}(\mathbf y)}{\partial y^{(j-1)}}+\frac{\partial\mu( \mathbf y)}{\partial y^{(j-1)}}F_{i}(\mathbf y )=\mu(\mathbf y )\frac{\partial F_j(\mathbf y)}{\partial y^{(i-1)}}+\frac{\partial \mu(\mathbf y )}{\partial y^{(i-1)}}F_{j}(\mathbf y),\; i=2,\cdots,n,\,j=1,2,\cdots,i-1,
\end{array}
\right.
\end{equation}
where $\mathbf y=(t,y,y^\prime,\cdots,y^{(n-1)})$.
In the general case, to solve such system of partial differential equations is not an easy problem. Therefore, we consider some special cases of $\mu(t,y,y^\prime,\cdots,y^{(n-1)})$. In particular, we are look for an integrating factor $\mu( \xi)$, where $\xi:=\xi(t,y,y^\prime,\cdots,y^{(n-1)})=\alpha(t)\displaystyle\prod_{k=1}^n\alpha_k(y^{(k-1)})$. Here, we assume the functions $\alpha(t)$ and $\alpha_k\left(y^{(k-1)}\right), k=1,2,\cdots,n$ to be differentiable functions. By substituting $\mu(\xi)$ 
in  \eqref{sys1}, we get
\begin{equation}\label{sys2}
\left\{
\begin{array}{ll}
 \mu(\xi)\frac{\partial F_{i}(\mathbf y)}{\partial t}+\mu^\prime(\xi) \xi_t F_{i}(\mathbf y )=\mu(\xi )\frac{\partial F_0(\mathbf y)}{\partial y^{(i-1)}}+\mu^\prime(\xi)\xi_{ y^{(i-1)}}F_{0}(\mathbf y), i=1,2,\cdots,n,&\\\\
 \mu(\xi)\frac{\partial F_{i}(\mathbf y)}{\partial y^{(j-1)}}+\mu^\prime(\xi)\xi_ {y^{(j-1)}}F_{i}(\mathbf y )=\mu(\xi )\frac{\partial F_j(\mathbf y)}{\partial y^{(i-1)}}+ \mu^\prime(\xi) \xi_{ y^{(i-1)}}F_{j}(\mathbf y), i=2,\cdots,n,j=1,2,\cdots,i-1,
\end{array}
\right.
\end{equation}
where $\mu^\prime(\xi)=\frac{d\mu}{d\xi}$ and $\xi_\eta=\frac{\partial \xi}{\partial \eta}.$
Equivalently, we have 
\begin{equation}\label{sys22}
\left\{
\begin{array}{ll}
\displaystyle \frac{\mu^\prime(\xi)}{\mu(\xi)}=\frac{\frac{\partial F_0(\mathbf y)}{\partial y^{(i-1)}}-\frac{\partial F_{i}(\mathbf y)}{\partial t}}{\xi_t F_{i}(\mathbf y)-\xi_{y^{(i-1)}}F_{0}(\mathbf y)},\; i=1,2,\cdots,n,\\\\
\displaystyle \frac{\mu^\prime(\xi)}{\mu(\xi)}=\frac{\frac{\partial F_j(\mathbf y)}{\partial y^{(i-1)}}-\frac{\partial F_{i}(\mathbf y)}{\partial y^{(j-1)}}}{\xi_{y^{(j-1)}} F_{i}(\mathbf y)-\xi_{y^{(i-1)}}F_{j}(\mathbf y)},\; i=2,\cdots,n,\; j=1,2,\cdots,i-1.
\end{array}
\right.
\end{equation}
Hence, an integrating factor $\mu(\xi)$ of  \eqref{HODE} exists if 
\[
\displaystyle\frac{\frac{\partial F_0(\mathbf y)}{\partial y^{(i-1)}}-\frac{\partial F_{i}(\mathbf y)}{\partial t}}{\xi_t F_{i}(\mathbf y)-\xi_{y^{(i-1)}}F_{0}(\mathbf y)},\; i=1,2,\cdots,n
\]
and 
\[
\displaystyle \frac{\frac{\partial F_j(\mathbf y)}{\partial y^{(i-1)}}-\frac{\partial F_{i}(\mathbf y)}{\partial y^{(j-1)}}}{\xi_{y^{(j-1)}} F_{i}(\mathbf y)-\xi_{y^{(i-1)}}F_{j}(\mathbf y)},\; i=2,\cdots,n,\; j=1,2,\cdots,i-1
\]
are functions in $\xi$ and they are equal. Therefore, we have the following theorem:
%%%%%%%%%%%%%%%%%%%%%%%%%%%%%%%%Thm%%%%%%%%%%%%%%%%%%%%%%%%%%%%%%%%
\begin{theorem}\label{maintheorem} Let $\xi=\alpha(t)\displaystyle\prod_{k=1}^n\alpha_k(y^{(k-1)})$, where $\alpha(t)$ and $\alpha_k\left(y^{(k-1)}\right), k=1,2,\cdots,n$ are differentiable functions. Assume that $F_0,F_1,F_2,\ldots,F_n$ are continuous functions with their first partial derivatives on a simply connected domain $\Omega\subset \R^{n+1}$. Moreover, assume that Equation \eqref{HODE} is a none exact differential equation. Then it  admits a none constant integrating factor  $\mu(\xi)$   
if 
\[
\displaystyle\frac{\frac{\partial F_0(\mathbf y)}{\partial y^{(i-1)}}-\frac{\partial F_{i}(\mathbf y)}{\partial t}}{\xi_t F_{i}(\mathbf y)-\xi_{y^{(i-1)}}F_{0}(\mathbf y)},\; i=1,2,\cdots,n
\]
and 
\[
\displaystyle \frac{\frac{\partial F_j(\mathbf y)}{\partial y^{(i-1)}}-\frac{\partial F_{i}(\mathbf y)}{\partial y^{(j-1)}}}{\xi_{y^{(j-1)}} F_{i}(\mathbf y)-\xi_{y^{(i-1)}}F_{j}(\mathbf y)},\; i=2,\cdots,n,\; j=1,2,\cdots,i-1
\]
are equal and they are functions in $\xi$. Moreover, the integrating factor is  explicitly  given by 
\[
\mu(\xi)=\exp\left\{\displaystyle\int\frac{\frac{\partial F_0(\mathbf y)}{\partial y}-\frac{\partial F_{1}(\mathbf y)}{\partial t}}{\xi_t F_{1}(\mathbf y)-\xi_{y}F_{0}(\mathbf y)}d\xi\right\}.\; \square
\]
\end{theorem}
%%%%%%%%%%%%%%%%%%%%%%%%%%%%%%%%%%%%%END THM - COR%%%%%%%%%%%%%%%%%%%
\begin{corollary}  
Let $\xi=\alpha(t)$, where $\alpha(t)$ is a differentiable function in $t$. Assume that $F_0,F_1,F_2,\ldots,F_n$ are continuous functions with their first partial derivatives on a simply connected domain $\Omega\subset \R^{n+1}$. Moreover, assume that Equation \eqref{HODE} is a none exact differential equation. Then it  admits a none constant integrating factor  $\mu(\alpha(t))$   
if 
\begin{itemize}
\item [I)] for $i=2,\cdots,n$ and for $ j=1,2,\cdots,i-1$, we have
\[
\frac{\partial F_j(\mathbf y)}{\partial y^{(i-1)}}=\frac{\partial F_{i}(\mathbf y)}{\partial y^{(j-1)}},\; 
\]
and 
\item [II)] for $ i=1,2,\cdots,n$,  
\[
\displaystyle\left[\frac{\partial F_0(\mathbf y)}{\partial y^{(i-1)}}-\frac{\partial F_{i}(\mathbf y)}{\partial t}\right]/\left[\alpha^\prime(t) F_{i}(\mathbf y)\right],
\] 
are equal and they are functions in $\xi=\alpha(t)$. 
\end{itemize}
Moreover, the integrating factor is  explicitly  given by 
\[
\mu(\xi)=\exp\left\{\displaystyle\int\left[\frac{\partial F_0(\mathbf y)}{\partial y}-\frac{\partial F_{1}(\mathbf y)}{\partial t}\right]/\left[\alpha^\prime(t) F_{1}(\mathbf y)\right]d\xi\right\}.\; \square
\]
\end{corollary}
%%%%%%%%%%%%%%%%%%%%%%%%%%%%%%%%%%COR-EXAP%%%%%%%%%%%%%%%%%%%%%%%%%
\begin{example}
Consider the following $n$-th order linear differential equation:
\[
P_n(t) y^{(n)}+P_{n-1}(t)y^{(n-1)}+\cdots+P_2(t) y^{\prime\prime}+P_1(t)y^\prime+P_0(t)y=h(t),
\]
where $P_i(t),\;i=0,1,2,\cdots,n$  are non-zero differentiable functions on some open interval $(a,b)\subset \R,$ and $h(t)$ is continuous function on $(a,b).$
Then $F_n=P_n(t)$, $F_{(n-1)}=P_{(n-1)}(t),\cdots,F_1=P_1(t),\;F_0=P_0(t)y-h(t)$.
Clearly, 
\[
\frac{\partial F_j(\mathbf y)}{\partial y^{(i-1)}}=\frac{\partial F_{i}(\mathbf y)}{\partial y^{(j-1)}},\; i=2,\cdots,n\; \text{and}\; j=1,2,\cdots,i-1.
\]
Moreover, $\frac{\partial F_0(\mathbf y)}{\partial y^{(i-1)}}=0, \; \forall i=2,\cdots,n$, $\frac{\partial F_i(\mathbf y)}{\partial y}=P^\prime_i(t),\; \forall i=1,\cdots,n$, and $\frac{\partial F_0(\mathbf y)}{\partial y}=P_0(t).$ Hence, to have an integrating factor in $t$, we must have
\[
\frac{P^\prime_n(t)}{P_n(t)}=\cdots=\frac{P^\prime_2(t)}{P_2(t)}=\frac{P^\prime_1(t)-P_0(t)}{P_1(t)}.
\]
Therefore, $P_n(t),\cdots,P_2(t)$ must be linearly dependent. Moreover, $P_0$ and $P_1$ must satisfy  $W(P_1,P_2)(t)=P_0(t)P_2(t),$ where $W(P_1,P_2)$ is the Wronskian's of $P_1$ and $P_2$. 
In this case, the integrating factor $\mu(t)=\frac{1}{P_n(t)}.$ Hence, for non-zero and differentiable functions $P(t)$, $P_1(t)$ and $P_0(t)$ the differential equation
\[
a_nP(t) y^{(n)}+a_{n-1}P(t)y^{(n-1)}+\cdots+a_2 P(t) y^{\prime\prime}+P_1(t)y^\prime+P_0(t)y=h(t),
\]
has an integrating factor $\mu(t)=\frac{1}{P(t)}$ provided that $W(P_1,P)(t)=P_0(t)P(t).$
Using this condition, we have that $\frac{P_0}{P}=\left(\frac{P1}{P}\right)^\prime$. Therefore, the above differential equation becomes
\[
a_n y^{(n)}+a_{n-1}y^{(n-1)}+\cdots+ a_2 y^{\prime\prime}+\frac{P_1(t)}{P(t)}y^\prime+\left(\frac{P_1(t)}{P(t)}\right)^\prime y=\frac{P(t)}{h(t)},
\]
  and so, its first integral is given by
\[
a_n y^{(n-1)}+a_{n-1}y^{(n-2)}+\cdots+ a_3 y^{\prime\prime}+a_2y^\prime+\frac{P_1(t)}{P(t)}y=\int^t\frac{h(s)}{P(s)}ds+c,
\]
where $c$ is the integrating constant.
\end{example}
%%%%%%%%%%%%%%%%%%%%%%%%%%EXAP-COR%%%%%%%%%%%%%%%%%%%%%%%%%%%%%%%%%
\begin{corollary} For $k=1,2,\cdots,n$, let $\xi=\alpha_k(y^{(k-1)})$, where  $\alpha_k(y^{(k-1)})$ is a differentiable function in $y^{(k-1)}$. Assume that $F_0,F_1,F_2,\ldots,F_n$ are continuous functions with their first partial derivatives on a simply connected domain $\Omega\subset \R^{n+1}$. In addition, assume that Equation \eqref{HODE} is a none exact differential equation. Then it  admits a none constant integrating factor  $\mu(\alpha_k(y^{(k-1)}))$ if 
\begin{itemize}
\item [I)] $\displaystyle \frac{\partial F_i}{\partial t}=\frac{\partial F_0}{\partial y^{(i-1)}}$ for $i=1,2,\cdots,n$ and
 $\displaystyle\frac{\partial F_j(\mathbf y)}{\partial y^{(i-1)}}=\frac{\partial F_{i}(\mathbf y)}{\partial y^{(j-1)}},\;i=2,3,\cdots,n,\;j=1,2,\cdots, i-1,\; i,j\neq k; 
$
and 
\item [II)]    
\[
\displaystyle\left[\frac{\partial F_k(\mathbf y)}{\partial y^{(i-1)}}-\frac{\partial F_{i}(\mathbf y)}{\partial y^{(k-1)}}\right]/\left[\alpha_k^\prime(y^{(k-1)}) F_{i}(\mathbf y)\right],\; i=1,2,\cdots,n,\;i\neq k,
\]
and
\[
\displaystyle\left[\frac{\partial F_k(\mathbf y)}{\partial t}-\frac{\partial F_{0}(\mathbf y)}{\partial y^{(k-1)}}\right]/\left[\alpha_k^\prime(y^{(k-1)}) F_{0}(\mathbf y)\right]
\]
are equal and they are functions in $\xi=\alpha_k(y^{(k-1)})$. 
\end{itemize}
Moreover, the integrating factor is  explicitly  given by 
\[
\mu(\xi)=\exp\left\{\displaystyle\int\left[\frac{\partial F_k(\mathbf y)}{\partial t}-\frac{\partial F_{0}(\mathbf y)}{\partial y^{(k-1)}}\right]/\left[\alpha_k^\prime(y^{(k-1)}) F_{0}(\mathbf y)\right]d\xi\right\}.\;\square
\]
\end{corollary}
%%%%%%%%%%%%%%%%%%%%%%%%%%%%%%%%%%%%%End COR%%%%%%%%%%%%%%%%%%%%%%%%
\begin{example}
Consider the following third order differential equation:
\begin{equation}\label{ex3eq}
y^3y^{\prime\prime\prime}+y^3y^{\prime\prime}-2 t y^\prime +y=0.
\end{equation}
Then $F_3=y^3$, $F_2=y^3$, $F_1=-2t$, and $F_0=y.$ Hence, $F_{2y^{\prime\prime}}=F_{3y^\prime}=0$, $F_{0y^{\prime\prime}}=F_{3t}=0$,
$F_{0y^\prime}=F_{2t}=0$, $0=F_{1y^{\prime}}\neq F_{2y}=3y^2$, $0=F_{1y^{\prime\prime}}\neq F_{3y}=3y^2$, and $-2=F_{1t}\neq F_{0y}=1$. By applying the above corollary, an integrating factor in terms of $y$ exists for this equation. Particularly, $\mu(y)=y^{-3}$. By Multiplying \eqref{ex3eq} by $\mu(y)=y^{-3}$, we get
\[
y^{\prime\prime\prime}+y^{\prime\prime}-2 ty^{-3} y^\prime +y^{-2}=0.
\]
Clearly, this differential equation is exact. Moreover, its first integral is 
\[
y^{\prime\prime}+y^\prime+ty^{-2}=c,
\]
where $c$ is an integrating constant.
\end{example}%%%%%%%%%%%%%%%%%%%%%%%%%%%%%%%%%%%%%%%%%%%%%%%%%%%%%%%%%%%%%%%%%%%%%%
\section{Concluding Remarks}
In this paper, we investigated the existence of integrating factors of the following class of third order differential equations: 
\begin{align}\label{HODEcr}
F_n\left(t,y,y^\prime,y^{\prime\prime},\ldots,y^{(n-1)}\right)y^{(n)}&+F_{n-1}\left(t,y,y^\prime,y^{\prime\prime},\ldots,y^{(n-1)}\right)y^{(n-1)}+\cdots+ \nonumber\\
&+F_{1}\left(t,y,y^\prime,y^{\prime\prime},\ldots,y^{(n-1)}\right)y^{\prime}+F_{0}\left(t,y,y^\prime,y^{\prime\prime}\ldots,y^{(n-1)}\right)\nonumber\\
&=0.
\end{align}
where $F_0,F_1,F_2, \cdots,F_n$ are continuous functions with their first partial derivatives 
on some simply connected domain $\Omega \subset\R^{n+1}$. Particularly, we proved some results related to the existence 
 of  integrating factors of \eqref{HODEcr}. We also presented some illustrative examples. We remark that these results not only useful for finding integrating factors for \eqref{HODEcr} analytically but also computationally. In fact, we can check the validity of the conditions in our results by using the symbolic toolboxes in different mathematical softwares, e.g.,   MAPLE and MATLAB softwares. Also, by using these symbolic toolboxes, we can find the integrating factors of \eqref{HODEcr} by using the explicit forms  given in our results. We also remark that  by using the same approach in this paper, we can derive integrating factors of \eqref{HODEcr} in terms of $\xi=\alpha(t)+\beta(y)+\gamma(y^\prime)+\delta(y^{\prime\prime})$, where $\alpha(t)$, $\beta(t)$, $\gamma(t)$ and $\delta(t)$ are differentiable functions.
%_____________________________________________________________________------------------The Bibliography------------------------------------

%----------------------------------------------------------------------------------
\end{document}